\documentclass[reqno,10pt]{article}
\usepackage{hyperref}
\usepackage{amsmath}
\usepackage{amssymb}
\usepackage{mathrsfs}
\usepackage{multirow}

\newtheorem{thm}{Theorem}[section]
\newtheorem{lemma}[thm]{Lemma}

\newtheorem{corol}[thm]{Corollary}
\newtheorem{propos}[thm]{Proposition}
\newtheorem{rema}{Remark}[section]
\def\bp{\begin{propos}}
\def\ep{\end{propos}}
\def\bt{\begin{thm}}
\def\et{\end{thm}}
\def\bco{\begin{corol}}
\def\eco{\end{corol}}
\def\bl{\begin{lemma}}
\def\el{\end{lemma}}
\def\br{\begin{rema}}
\def\er{\end{rema}}
\def\be{\begin{equation}}
\def\ee{\end{equation}}
\def\ba{\begin{array}}
\def\ea{\end{array}}
\def\bena{\begin{eqnarray}}
\def\eena{\end{eqnarray}}

\def\P{{\mathbb P}}
\def\E{{\mathbb E}}
\def\R{{\mathbb R}}
\def\Z{{\mathbb Z}}

\def\1{I}
\def\imath{\textbf{i}}
\def\jmath{\textbf{j}}
\def\chi{\zeta}

\def\a{{\alpha}}

\def\Var{\hb{Var}}

\def\QED{\hfill$\square$\vskip 3mm}

\def\hb{\hbox}
\def\({\left(}
\def\){\right)}

\begin{document}

\title{On The Time Constant for Last Passage Percolation on Complete Graph
\\[5mm]
\footnotetext{*Correspondence author}
\footnotetext{AMS classification (2010): 60K 35} 
\footnotetext{Research supported in part by the Natural Science Foundation of China (under
grants 11271356, 11471222, 11671275)}}

\author{Xian-Yuan Wu$^*$,\ \ \ Rui Zhu}\vskip 10mm
\date{}
\maketitle

\begin{center}
\begin{minipage}{13cm}

\noindent School of Mathematical Sciences, Capital
Normal University, Beijing, 100048, China. Email:
\texttt{wuxy@cnu.edu.cn},\ \texttt{1073755862@qq.com}
\end{minipage}
\end{center}
\vskip 5mm
{\begin{center} \begin{minipage}{13cm}
{\bf Abstract}: This paper focuses on the time constant for last passage percolation on complete graph. Let $G_n=([n],E_n)$ be the complete graph on vertex set $[n]=\{1,2,\ldots,n\}$, and i.i.d. sequence $\{X_e:e\in E_n\}$ be the passage times of edges. Denote by $W_n$ the largest passage time among all self-avoiding paths from 1 to $n$. First, it is proved that $W_n/n$ converges to constant $\mu$, where $\mu$ is called the {\it time constant} and coincides with the {\it essential supremum} of $X_e$. Second, when $\mu<\infty$, it is proved that the deviation probability $\P(W_n/n\leq \mu-x)$ decays as fast as $e^{-\Theta(n^2)}$, and as a corollary, an upper bound for the variance of $W_n$ is obtained. Finally, when $\mu=\infty$, lower and upper bounds for $W_n/n$ are given.


\end{minipage}
\end{center}}

 \vskip 5mm
\section{Introduction and statement of the results}
\renewcommand{\theequation}{1.\arabic{equation}}
\setcounter{equation}{0}

The last passage percolation process has been widely studied over the last two decades (see for example \cite{BDMMZ,BS,BM,HM,J,M} and the references therein). There are several equivalent physical interpretations for the model (see \cite{GW,O',W}), this makes the model playing important roles in statistical physics. With interest arisen from the study of real-world networks, in this paper, we will consider the last passage percolation model on complete graph. For a network with links carrying an uncertain cost, we try to characterize the {\it worst} path within it. In other words, when the link weight is interpreted as some kind of income, the goal is just to characterize the optimal path in the network. This paper is motivated by the three resent papers for first passage percolation on complete graph from Eckhoff, Goodman, Hofstad and Nardi \cite{EGHN1,EGHN2,EGHN3}, actually, we are running on the same way, but on the opposite direction.

We start by introducing the model. Let $G_n=([n],E_n)$ be the complete graph on vertex set $[n]=\{1,2,\ldots,n\}$, where edge set $E_n=\{\langle i,j\rangle:1\leq i<j\leq n\}$. To each edge $e\in E_n$ assign a random edge passage time $X_e$, where $\{X_e:e\in E_n\}$ are positive, independent and identically distributed random variables. The set of all self-avoiding paths between vertex 1 and n is denoted by $\Pi_{1,n}$. The passage time of a path $\pi$ in $\Pi_{1,n}$ is a function $T(\pi)$ defined as follows $$T(\pi)=\sum_{e\in\pi}X_e.$$ Let $W_n$ be the largest passage time among all self-avoiding paths from 1 to n, i.e., \be\label{1}W_n=\sup_{\pi\in\Pi_{1,n}}T(\pi).\ee

To explore the asymptotic behavior of $W_n$ (if necessary, properly scaled) for large $n$ is the main goal of the present paper. Clearly, in this point of view, the last passage percolation model on complete graph is mathematically equivalent to the first passage model, provided {\it negative} edge passage times are allowed.

For first passage percolation model on $G_n$ with non-negative edge passage time (distribution not scaled in $n$), it is easy to see that the corresponding {\it time constant} always exists, but equals zero, this indicates that the time constant problem is no longer a problem. But the situation for last passage percolation is different, first of all, we have the following existence theorem for time constant.
\bt\label{t1} For any distribution of edge passage time, the time constant of the model exists. More precisely, there exists some constant $0<\mu\leq \infty$, such that \be\label{3}
\frac{W_n}n\rightarrow\mu\ \ \rm {a.s.}\ee as $n\rightarrow \infty$. In particular, when $\mu<\infty$, the above convergence is also in $L_1$. Furthermore, $\mu$ coincides with the essential supremum of $X_e$, i.e.
\be\label{2}\mu=\inf\{x:\P(X_e>x)=0\}.\ee  \et

Now, let $F(x)=\P(X_e\leq x)$ and $H(x)=1-F(x)$, $x\in\R$ be the distribution function and the tail probability function of $X_e$. By the above Theorem~\ref{t1}, when $\mu=\infty$, $W_n/n$ tends to $\infty$ as $n\rightarrow \infty$. Here, we give lower and upper bounds to $W_n/n$ as in the following theorem.

\bt\label{t2}Suppose that $\mu=\infty$. Let $f$ and $g$ be two functions such that $H(f(n))=\ln n/n$ and $n^2H(g(n))\rightarrow 0$ as $n\rightarrow\infty$. Then, for any $\epsilon>0$
\be\lim_{n\rightarrow\infty}\P\((1-\epsilon)f(n)\leq \frac{W_n}n\leq g(n)\)=1.\ee \et

Obviously, we have the following corollary for exponential and power-law tail probability function $H$. We don't think the bounds given in Theorem~\ref{t2} are sharp, but i) of Corollary~\ref{co} indicates that they are not so bad.
\bco\label{co}Suppose that $H(x)=\P(X_e>x)$, $x\in \R$. Then
\begin{description}
  \item[\ \ \ \ i)] if $H(x)\sim e^{-\lambda x}$, $\lambda>0$, then for any $\epsilon>0$
  \be\lim_{n\rightarrow\infty}\P\((1-\epsilon)\frac{\ln n-\ln\ln n}{\lambda}\leq \frac{W_n}n\leq \frac{2\ln n+\ln\ln n}{\lambda}\)=1;\ee
  \item[\ \ \ \ ii)] if $H(x)\sim x^{-\a}$, $\a>0$, then
  \be\lim_{n\rightarrow\infty}\P\((1-\epsilon)n^{\frac 1\alpha}[\ln n]^{-\frac 1\alpha}\leq \frac{W_n}n\leq {n^{\frac 2\alpha}}[\ln\ln n]^{\frac 1\alpha}\)=1.\ee

\end{description}

\eco

Finally, in the case where $\mu<\infty$, we are interested in the asymptotic behavior of the deviation probability $\P(W_n/n\leq \mu-x)$, $0<x<\mu$, for large $n$. In fact, we have

\bt\label{t3}
Suppose $\mu<\infty$ and $H$ be the tail probability of $X_e$. Then, for any $0<x<\mu$ with $p=p(x)=H(\mu-x)<1$, there exists some constants $C_1,\ c_1>1$ and $c_2>C_2>0$ such that
\be\label{4}c_1^ne^{-c_2n^2}\leq\P\(\frac {W_n}n\leq \mu-x\)\leq C_1^n e^{-C_2n^2}\ee for large enough $n$.
\et

In the case where $\mu<\infty$, Theorem~\ref{t3} has the following corollary for the variance of $W_n$. Note that, for first passage percolation on $\Z^d,\ d\geq 2$, the corresponding problem has been closely followed with interests by mathematicians and physicists for a long period of time. For details on this aspect, please refer to \cite{BeKS,DHS,D}.

\bco\label{co'} Suppose that $\mu<\infty$. Then we have the following upper bounds for the variance of ${W_n}$:
\begin{description}
  \item[\ \ \ i)] in the case when $\P(X_e=\mu)=p_0>0$, there exists some constant $D_1>0$, such that
  \be\label{c1}\Var\({W_n}\)\leq D_1n{\ln n}. \ee
  \item[\ \ \ ii)] in the case when the tail probability function $H$ is left continuous at $\mu$, let $\bar x(n)$ be the unique solution of equation ${\ln n}/n=xH(\mu-x/2)$, then there exists some constant $D_2>0$ such that \be\label{c2}\Var\({W_n}\)\leq D_2 \bar x(n)n^2.\ee For example, if $H(y)=(\mu-y)^\alpha$ for large $y<\mu$, $\alpha>0$, then we have
      $$ \Var\({W_n}\)\leq D_2 \bar x(n)n^2\leq 2^{\frac{\alpha}{\alpha+1}}D_2 n^{2-\frac 1{\alpha +1}}\ln^{\frac 1{\alpha+1} }n. $$
\end{description}

\eco

\br It is believed that the variance of $W_n$ is at least sublinear, but the authors can not give a proof to this declaration at the present time.\er

The rest of the paper is arranged as follows. In Section 2, first, by using the subadditive ergodic theorem of Liggett's version \cite{L}, we prove the existence of the time constant; second, by using a famous theorem of Ajtai, Koml\'os and Szemer\'edi \cite{AKS} on long path in random graph (see also \cite[Theorem 8.1]{B}), we characterize the time constant and then give a proof to Theorem~\ref{t2}. In Section 3, we first introduce the well known {\it Depth First Search} algorithm (see \cite{CLRS,KT}) for finding long paths in random graph, then, by using the constructed DFS random graph process, we prove Theorem~\ref{t3} and Corollary~\ref{co'}.

\vskip 5mm
\section{Proofs of Theorems~\ref{t1}, \ref{t2}.}
\renewcommand{\theequation}{2.\arabic{equation}}
\setcounter{equation}{0}
 In this section, we prove Theorem~\ref{t1} and Theorem~\ref{t2}. For any integers $0\leq m<n$, let $G_{m,n}$ be the complete graph on vertex set $\{m+1,m+2,\ldots,n\}$. Consider the last passage percolation model on $G_{m,n}$ and denote by $W_{m,n}$ the largest passage time from $m+1$ to $n$. Then, one has $G_n=G_{0,n}$ and $W_n=W_{0,n}$.

 {\it Proof of Theorem~\ref{t1}, the existence part:} Let us consider the two index process $\{W_{m,n}:0\leq m<n\}$, by the definition of the last passage percolation, one has
 \begin{description}
 \item[\ \ \  1)] $W_{0,n}\geq W_{0,m}+W_{m,n}$, wherever $0<m<n$;
  \item[\ \ \  2)] The joint distributions of $\{W_{m+1,m+k+1}:k\geq 1\}$ are the same as those of $\{W_{m,m+k}:k\geq 1\}$ for each $m\geq 0$;
  \item[\ \ \  3)] For each $k\geq 1$, $\{W_{nk,(n+1)k}:n\geq 1\}$ is an i.i.d. random variable sequence and hence an ergodic process.
 \end{description}
 By items 1) and 2) above and the classical subadditive limit theorem, one has $$\lim_{n\rightarrow\infty}\frac {\E(W_{0,n})}n=\sup_n\frac {\E(W_{0,n})}n=\mu\leq \infty.$$
 Then, in the case where $\mu<\infty$, the two index process $\{W_{m,n}:0\leq m<n\}$ satisfies the following property additionally

 {\bf 4}) For each $n$, $\E(W_{0,n})<\infty$ and $\E(W_{0,n})\leq cn$ for some constant $c(=\mu)$.

 \noindent Clearly,  $\{W_{m,n}:0\leq m<n\}$ always satisfies the following property

 {\bf 4'}) $\E(W_{0,1}^-)<\infty$.

 By items 1)-4) above and the subadditive ergodic Theorem of Liggett's version (see \cite[Theorem 1.10]{L}), one has $W_n/n\rightarrow\mu$ a.s. and in $L_1$ as $n\rightarrow\infty$.

 In the case where $\mu=\infty$, one has no $L_1$ convergence. The corresponding $a.s.$ convergence can be obtained as an extension of Theorem~1.10 in \cite{L} by a simple truncation argument as in the proof of Theorem 1.8 in \cite{K}, with 4) replaced by 4').\QED

 To finish the remain part of Theorem~\ref{t1}, we will first introduce the famous existence theorem for long path in a super-critical random graph due to Ajtai, Koml\'os and Szemer\'edi \cite{AKS}. Denote by $G_{n,p}$ the random subgraph of the complete graph $G_n$ with vertex set $[n]$ and random edge set ${\cal E}=\{e\in E_n:Y_e=1\}$, where $\{Y_e:e\in E_n\}$ is the i.i.d. Bernoulli random variable sequence with parameter $p$. In 1960, Paul Erd\H{o}s and Alfr\'ed R\'enyi \cite{ER} made the following
fundamental discovery: the random graph $G_{n,p}$ undergoes a phase transition around the edge probability $p=p(n) =1/n$. For any constant $\epsilon> 0$, if $p={(1-\epsilon)}/n$, then, $G_{n,p}$ has {\bf whp} all connected components of size at most logarithmic in $n$, while for $p={(1+\epsilon)}/ n$ {\bf whp} a unique connected component of linear size emerges in $G_{n,p}$. Where {\bf whp} means {\it with probability tends to 1 as $n\rightarrow\infty$}.

Although for the super-critical case $p={(1+\epsilon)}/ n$ the result of Erd\H{o}s and R\'enyi shows a typical
existence of a linear sized connected component, it does not imply that a longest path in such a
random graph is {\bf whp} linearly long. This was established about 20 years later by Ajtai, Koml\'os and
Szemer\'edi \cite{AKS}, the following is a version of this result which can be found in \cite[Theorem 8.1]{B}.

\bt\label{t2.1}\cite{B} Let $0<\theta=\theta(n)<\ln n-3\ln\ln n$ and $p=\frac \theta n$. Then $G_{n,p}$ contains {\bf whp} a path of length at least $$\(1-\frac {4\ln2}\theta \)n.$$\et

 By using Theorem~\ref{t2.1}, we finish the proof of Theorem~\ref{t1} as follows.

 {\it Proof of Theorem~\ref{t1}, the remain part:} Denote by $U^{ess}_X$ the essential supremum of $X_e$. In the case where $U^{ess}_X=\infty$, for any given $M>0$, let $p=p(M)=H(2M)>{\ln n}/{2n}$. By Theorem~\ref{t2.1} with $\theta(n)={\ln n}/ 2$, $G_{n,p}$ contains {\bf whp} a path of length at least $\(1-{4\ln2}/\theta \)n.$ This means that, for last passage percolation on $G_{n}$, there exists {\bf whp} a path, say $\pi$, of length $\(1- {4\ln2}/\theta \)n$, such that for each $e\in \pi$, $X_e>2M$.

 Based on the long path $\pi$, now we go to construct a path $\bar\pi\in\Pi_{1,n}$, such that $\bar\pi$ differs from $\pi$ for only several edges. Denote by ${\cal E}(1)$, ${\cal E}(n)$ the set of edges with $1$, $n$ as one of its endvertex respectively. Suppose the path $\pi$ is divided into $i_0$ pieces, say $\pi_1,\pi_2,\ldots,\pi_{i_0}$, after $1,n$ and all edges in ${\cal E}(1)$, ${\cal E}(n)$ are deleted. Clearly, one has $1\leq i_0\leq 3$. Denote by $u_j,v_j$ the endpoints of path $\pi_j$, $1\leq j\leq i_0$. Let $\bar\pi$ be the resulted self-avoiding path after $1,n$ and $\pi_1,\pi_2,\ldots,\pi_{i_0}$ are connected by edges $\langle 1,u_1\rangle, \langle v_1,u_2\rangle,\ldots,\langle v_{i_0-1},u_{i_0}\rangle$ and $\langle v_{i_0}, n\rangle$.


 By the construction of $\bar\pi$, there are at least $\(1-{4\ln2}/\theta \)n-4$ common edges between $\bar\pi$ and $\pi$, then $T(\bar\pi)$, the passage time of $\bar\pi$, $\geq [\(1-{4\ln2}/\theta \)n-4]\times 2M>Mn.$ So, {\bf whp}, one has
 $$\frac{W_n}n>M.$$
By the existence part of Theorem~\ref{t1}, one has $\mu>M$ and then $\mu=U^{ess}_X=\infty$.

In the case where $U^{ess}_X<\infty$, for any given $\delta>0$, let $p=p(\delta)=H(U^{ess}_X-\delta)>{\ln n}/{2n}=\theta/ n$. By using the same argument as above, {\bf whp}, one can construct a long path $\bar\pi\in\Pi_{1,n}$ such that $T(\bar\pi)\geq [\(1-{4\ln2}/\theta \)n-4]\times(U_{ess}-\delta)>(U_{ess}-2\delta)n$. So, one has {\bf whp}
 $$\frac{W_n}n>U^{ess}_X-2\delta.$$
By the existence part of Theorem~\ref{t1}, one has $\mu>U^{ess}_X-2\delta$ and then $\mu=U^{ess}_X$. \QED

By using the result of Theorem~\ref{t2.1}, we can also give a proof to Theorem~\ref{t2}.

{\it Proof of Theorem~\ref{t2}}: Suppose that $f(n)$ satisfy $H(f(n))={\ln n}/{n}$. By using Theorem~\ref{t2.1} to random graph $G_{n,p/2}$ and the monotonicity, we have in $G_{n,p}$, {\bf whp}, there exists a path $\pi$ of length at least $\(1-{8\ln 2}/{\ln n}\)n$. In other words, for last passage percolation on $G_n$, {\bf whp}, there exists a path $\pi$ of length at least $\(1-{8\ln 2}/{\ln n}\)n$, such that $X_e>f(n)$ for any $e\in \pi$. Construct a path $\bar \pi\in\Pi_{1,n}$ as done in the proof of Theorem~\ref{t1}, such that $T(\bar\pi)\geq [\(1-{8\ln 2}/{\ln n}\)n-4]f(n)$. Then, for any $\epsilon>0$, {\bf whp} $W_n\geq T(\bar\pi)\geq (1-\epsilon)f(n)n$. Thus we get the lower bound part of the theorem.

The upper bound part of the theorem is quite simple, and it follows from the fact that, if all edge passage times are less than $g(n)$, then $W_n/n\leq g(n)$. Actually, by the condition on $g(n)$, one has
$$\P(\exists\ e\ \rm {such\ that}\ X_e>g(n))\leq \binom{n}{2}H(g(n))\leq n^2H(g(n))\rightarrow 0$$ as $n\rightarrow\infty$.\QED

\vskip 5mm
\section{Depth First Search and the proof of Theorem~\ref{t3}.}
\renewcommand{\theequation}{3.\arabic{equation}}
\setcounter{equation}{0}

In this section we will prove Theorem~\ref{t3}.  In Section 2, we tried to construct a long bad {\it enough} path with high probability. But to finish the proof of Theorem~\ref{t3}, we are asked to construct a long {\it not so }bad path with probability tends to 1 fast {\it enough}. To this end, based on \cite{Kr}, we first introduce the {\it Depth First Search} algorithm for finding long path in graphs.


Let $G=([n],E)$ be a graph on vertex set $[n]=\{1,2,\ldots,n\}$. The algorithm receives as an input $G$, and maintains three sets of vertices, and a set of edges. Let $S$ be the set of vertices whose exploration is complete, $T$ be the set of unvisited vertices, and $U=[n]\setminus(S\cup T)$. Let $\hat E$ be the set of edges in $E_n$ whose exploration is complete, recall that $E_n$ is the edge set of complete graph $G_n$.

The algorithm initializes with $S_0=U_0=\emptyset,\ T_0=[n]$ and ${\hat E}_0=\emptyset$. Suppose that $(S_t,U_t,T_t,{\hat E}_t)$ be well defined, we define $(S_{t+1},U_{t+1},T_{t+1},{\hat E}_{t+1})$ as the following.

In the case when $U_t\not=\emptyset$, let $x$ be the last vertex that has been added to $U_t$, then we try to look for a vertex $y$ such that edge $\langle x,y\rangle\notin {\hat E}_t$. If such a vertex exists, let $y$ be the smallest one, then let $S_{t+1}=S_t$, $U_{t+1}=U_t\cup\{y\}$, $T_{k+1}=T_t\setminus\{y\}$ and ${\hat E}_{t+1}={\hat E}_t\cup\{\langle x,y\rangle\}$, if $y\in T_t$ and $\langle x,y\rangle\in E$; let
$S_{t+1}=S_t$, $U_{t+1}=U_t$, $T_{k+1}=T_t$ and ${\hat E}_{t+1}={\hat E}_t\cup\{\langle x,y\rangle\}$ otherwise. In the case when such a vertex does not exist, let $S_{t+1}=S_t\cup\{x\}$, $U_{t+1}=U_t\setminus\{x\}$, $T_{k+1}=T_t$ and ${\hat E}_{t+1}={\hat E}_t$.

In the case when $U_t=\emptyset,\ T_t\not=\emptyset$, denote by $x$ the smallest one in $T_t$ and then let $S_{t+1}=S_t$, $U_{t+1}=\{x\}$, $T_{k+1}=T_t\setminus\{x\}$ and ${\hat E}_{t+1}={\hat E}_t$.

The algorithm stops whenever $U_t=T_t=\emptyset$.

Observe that the DFS algorithm starts revealing a connected component $C$ of $G$ at the moment the smallest vertex of $C$ gets into $U$ and completes discovering all of $C$ when $U$ becomes empty again. The period of time between two consecutive emptying of $U$ is called an {\it epoch}, and each epoch corresponds to one connected component of $G$. Denote by $N_c(G)$ the number of connected components of $G$.

For the DFS algorithm defined above, the following properties are immediate to verify:
\begin{description}
  \item[\ \ \ \ (P1)] at each round of the algorithm, at most one edge of $E_n$ is checked (whether or not lies in $E$) and then enters $\hat E$. When no edge is checked, one vertex moves: either some vertex move from $U$ to $S$, or some vertex of $T$, which is the smallest one of some connected component of $G$, moves from $T$ to $U$. Then the DFS process stops at time $N=n+N_c(G)+\binom {n}{2}$.

  \item[\ \ \ \ (P2)] at any time $t$, each edge in $\{\langle u,v\rangle:u\in S_t, v\in T_t\}\subset E_n$ has been checked before $t$ to be not in $E$. This implies that $|{\hat E}_t|\geq |S_t|\times|T_t|$.
  \item[\ \ \ \ (P3)] the set $U$ always spans a path. This implies that, $G$ contains a path of length $$\sup_{1\leq t\leq N}|U_t|-1.$$
\end{description}

{\small

\begin{center}
\begin{minipage}{13cm}

\begin{tabular}{|c|c|c|c|c|}
        \hline
       Step &      S       &   U & T & $\hat E$ \\
       \hline
       0    & $\emptyset$ &$\emptyset$&$\{1,2,3,4,5\}$&$\emptyset$ \\
       \hline
        1    & $\emptyset$ &$\{1\}$&$\{2,3,4,5\}$&$\emptyset$ \\
        \hline
         2  & $\emptyset$ &$\{1\}$&$\{2,3,4,5\}$&$\{e_1=\langle 1,2\rangle\}$  \\
        \hline
        3  & $\emptyset$ &$\{1\}$&$\{2,3,4,5\}$&$\{e_1,e_2={\langle 1,3\rangle}\}$  \\
        \hline
        4  & $\emptyset$ &$\{1,4\}$&$\{2,3,5\}$&$\{e_1,e_2,e_3=\underline{\langle 1,4\rangle}\}$  \\
        \hline
        5  & $\emptyset$ &$\{1,4\}$&$\{2,3,5\}$&$\{e_1,e_2,e_3,e_6=\langle 2,4\rangle\}$  \\
        \hline
       6  & $\emptyset$ &$\{1,4\}$&$\{2,3,5\}$&$\{e_1,e_2,e_3,e_6,e_8=\langle 3,4\rangle\}$  \\
        \hline
        7  & $\emptyset$ &$\{1,4\}$&$\{2,3,5\}$&$\{e_1,e_2,e_3,e_6,e_8,e_{10}=\langle 4,5\rangle\}$  \\
         \hline
         8 & $\{4\}$ &$\{1\}$&$\{2,3,5\}$&$\{e_1,e_2,e_3,e_6,e_8,e_{10}\}$  \\
        \hline
        9 & $\{4\}$ &$\{1\}$&$\{2,3,5\}$&$\{e_1,e_2,e_3,e_6,e_8,e_{10},e_4=\langle 1,5\rangle\}$  \\
        \hline
        10 & $\{1,4\}$ &$\emptyset$&$\{2,3,5\}$&$\{e_1,e_2,e_3,e_6,e_8,e_{10},e_4\}$  \\
        \hline
         11 & $\{1,4\}$ &$\{2\}$&$\{3,5\}$&$\{e_1,e_2,e_3,e_6,e_8,e_{10},e_4\}$  \\
        \hline
         12 & $\{1,4\}$ &$\{2,3\}$&$\{5\}$&$\{e_1,e_2,e_3,e_6,e_8,e_{10},e_4,e_5=\underline{\langle 2,3\rangle}\}$  \\
        \hline
        13 & $\{1,4\}$ &$\{2,3,5\}$&$\emptyset$&$\{e_1,e_2,e_3,e_6,e_8,e_{10},e_4,e_5,e_9=\underline{\langle 3,5\rangle}\}$  \\
        \hline
         14 & $\{1,4\}$ &$\{2,3,5\}$&$\emptyset$&$\{e_1,e_2,e_3,e_6,e_8,e_{10},e_4,e_5,e_9,e_7=\underline{\langle2,5\rangle}\}$  \\
        \hline
        15 & $\{1,4,5\}$ &$\{2,3\}$&$\emptyset$&$\{e_1,e_2,e_3,e_6,e_8,e_{10},e_4,e_5,e_9,e_7\}$  \\
        \hline
         16 & $\{1,3,4,5\}$ &$\{2\}$&$\emptyset$&$\{e_1,e_2,e_3,e_6,e_8,e_{10},e_4,e_5,e_9,e_7\}$  \\
        \hline
         17 & $\{1,2,3,4,5\}$ &$\emptyset$&$\emptyset$&$\{e_1,e_2,e_3,e_6,e_8,e_{10},e_4,e_5,e_9,e_7\}$  \\
        \hline

 \end{tabular}

\vskip 3mm
 Figure 1: {\it Let $e_1,e_2,\ldots,e_9,e_{10}$ denote edge
 $\langle 1,2\rangle,\langle 1,3\rangle\ldots,\langle 3,5\rangle,\langle 4,5\rangle$ respectively. The above is the DFS process for graph $G=([5], E)$ with $E=\{\langle 1,4\rangle,\langle 2,3\rangle,\langle 2,5\rangle, \langle 3,5\rangle\}$.
 }

 \end{minipage}

 \end{center}
}
\vskip 5mm

The following proposition exploit the features of DFS algorithm to derive the existence of long path in graph $G$.

\bp\label{p3.1}\cite{BKS} Let $k<n$ be  positive integers. Assume that $G=([n],E)$ is a graph on $[n]$, containing an edge between any two disjoint subsets $S,T\subset [n]$ of size $|S|$=$|T|$=k. Then, $G$ contains a path of length $n-2k+1$.\ep
{\it Proof}: Run the DFS algorithm on $G$, let $\tau=\inf\{t:|S_t|=|T_t|\}$. Note that by property {\bf (P1)}, $\tau$ is well defined. Since $G$ has no edge between $S_{\tau}$ and $T_{\tau}$ by property {\bf (P2)}, it follows from the assumption of  the proposition that $|S_\tau|=|T_\tau|\leq k-1$. Then, by property {\bf(P3)}, $G$ contains a path of length at least $|U_\tau|-1=n-|S_\tau|-|T_\tau|-1\geq n-2k+1$.\QED

Now, we run the DFS algorithm on random graph $G_{n,p}=([n],\cal E)$ and define the random DFS process $\{(S_t,U_t,T_t,{\hat E}_t):0\leq t\leq N\}$, where $N=n+\binom{n}{2}+N_c(G_{n,p})$ is a random variable. Note that in the random case, at each round of the algorithm, when a edge $e\in E_n$ is checked, the content is changed from ``whether or not $e\in E$" to ``whether or not $e\in \cal E$".

{\it Proof of Theorem~\ref{t3}}: Let's consider the last passage percolation on $G_n$, assume that the time constant $\mu<\infty$. For any $x>0$, such that $p:=H(\mu-x)<1$, take $\epsilon>0$ small enough such that $x'=\mu-(\mu-x)/(1-2\epsilon)>0$. Let $p'=H(\mu-x')$, and consider the random graph $G_{n,p'}$.

Let $k=\lfloor\frac{\epsilon n}2\rfloor$, then the probability that there exists disjoint $S,T\subset [n]$, $|S|=|T|=k$, such that $G_{n,p'}$ contains no edge between $S,T$ is at most
\be\label{5}\binom{n}{k}\binom{n-k}k (1-p')^{k^2}\leq \binom nk^2e^{-p'k^2}\leq \(\frac {en}k\)^{2k}e^{-p'k^2}\leq C_1^n e^{-C_2 n^2},\ee where $C_1=\({2e}/\epsilon\)^\epsilon>1$ and $C_2={\epsilon^{2}}H(\mu-x')/{5}>0$.
On the other hand, if in $G_{n,p'}$ there exists long path with length at least $n-2k+1$, then as we have done in Section 2, this means that, for last passage percolation on $G_n$, $W_n\geq (n-2k-3)(\mu-x')$. Hence, by Proposition~\ref{p3.1} and (\ref{5}),
$$\ba{ll}\P(W_n\leq (\mu-x)n)&=\P(W_n\leq (1-2\epsilon)(\mu-x')n)\\
&\leq\P(W_n<(n-2k-3)(\mu-x'))\\
&\leq  C_1^n e^{-C_2 n^2}.
\ea
$$
We get the upper bound part of the theorem.

It is easy to see that, if the passage times of all edges are less than $\mu-x$, then $W_n\leq (\mu-x)n$. This implies that
$$\P(W_n\leq (\mu-x)n)\geq (1-p)^{\frac {n^2-n}2}=(1-p)^{-n/2}e^{-[\frac 12\ln \frac 1{1-p}] n^2}=c_1^ne^{-c_2 n^2}.$$
Thus, we get the lower bound and finish the proof of the theorem.\QED

{\it Proof of Corollary~\ref{co'}}: First of at all, by the upper bound given in Theorem~\ref{t3}, one has for small enough $x>0$,
\be\label{6}\ba{ll}\Var\({\frac {W_n}{n}}\)&=\E\(\frac {W_n}{n}\)^2-\E^2\(\frac {W_n}{n}\)\leq \mu^2-(\mu-x)^2\(1-C_1^ne^{-C_2 n^2}\)^2\\[3mm]
&\leq \mu^2-(\mu^2-2x)\(1-2C_1^ne^{-C_2 n^2}\)\\[3mm]
&\leq 2x+2\mu^2C_1^ne^{-C_2n^2},\ea\ee where $C_1=\({2e}/\epsilon\)^\epsilon$ and $C_2={\epsilon^{2}}H(\mu-x')/{5}$ are given in (\ref{5}).

In the case of $\P(X_e=\mu)=p_0>0$, one has $H(\mu-x)\geq p_0>0$ for all $0<x<\mu$.
Now, for small enough $x>0$, take $x'=x/2$, then $x/4\mu\leq\epsilon=x/(4\mu-2x)\leq x/3\mu$, and
\be\label{8} C_1\leq\(\frac{6e\mu}x\)^{x/{3\mu}},\ C_2\geq\frac {p_0}{80\mu^2}x^2.\ee
Note that here we have used the fact that $\(1/x\)^x$ is increasing in $x$, while $x<e$. Take $x=D{\ln n}/n$ with $D=27\mu/p_0$, then
$$C_1^ne^{-C_2n^2}\leq \exp\left\{-\(\frac {27^2}{80}-9\)\frac1{p_0}\ln^2 n\right\}.$$
Hence, by (\ref{6}),
\be\label{7}\Var\(\frac{W_n}n\)\leq D_1\frac{\ln n} n,\ee for some $D_1>2D$.

In the case when $H$ is left continuous at $\mu$, i.e. $\lim_{y\rightarrow \mu^-}H(y)=H(\mu)=0$. For small enough $x>0$, let $x'=x/2$, then $C_1$ has the same upper bound as in (\ref{8}) and $C_2\geq  {H(\mu-x/2)x^2}/{80\mu^2}$. Noticing that $xH(\mu-x/2)$ is strictly monotonic in $x$, we define $\bar x(n)$ be the unique solution of equation $\ln n/n=xH(\mu-x/2)$. Clearly $\bar x(n)\geq \ln n/n.$

Now, Let $x=C_3\bar x(n) $ with $C_3=27\mu$, then $C_4:=C_3^2/80\mu^2-C_3/3\mu>0$ and
$$ C_1^ne^{-C_2n^2}\leq e^{-C_4\bar x(n)n\ln n}\leq n^{-C_4\ln n}.$$
Hence, by (\ref{6}),
\be\label{9}\Var\(\frac{W_n}n\)\leq D_2\bar x(n),\ee for some $D_2>2C_3$.\QED


\vskip10mm

\end{document}